\documentclass{elsart}

\usepackage{amssymb}
\makeatletter\@addtoreset{equation}{section}\makeatother
\def\theequation{\arabic{section}.\arabic{equation}}
\newtheorem{Tm}{Theorem}[section]
\newtheorem{Rk}{Remark}[section]
\newtheorem{Lm}{Lemma}[section]
\newtheorem{Co}{Corollary}[section]
\newtheorem{Df}{Definition}[section]
\newtheorem{Pn}{Proposition}[section]

\begin{document}

\begin{frontmatter}

\title{Infinite Order Differential Operators in Spaces of Entire
Functions\thanksref{titlefn}}
\thanks[titlefn]{Supported in part by KBN grant 2P03A 02915}

\author[Lublin]{Yuri Kozitsky\thanksref{and},}
\ead{jkozi@golem.umcs.lublin.pl kozitsky@physik.uni-bielefeld.de }
\author[Lublin]{Piotr Oleszczuk}
\ead{poleszcz@golem.umcs.lublin.pl}
\author[Davis]{and Lech Wo{\l}owski}
\thanks[and]{Address for the period October 2002 - June 2003: Research Center
 BiBoS, Universit{\"a}t Bielefeld, Postfach 100131, D 33501 Bielefeld, Germany  }
\ead{wolowski@math.ucdavis.edu}
\address[Lublin]{Instytut Matematyki,
Uniwersytet Marii Curie-Sk{\l}odowskiej, PL 20-031 Lublin, Poland}
\address[Davis]{Department of Mathematics, University of California, Davis, CA 95616-8633, USA}
\begin{abstract}

Differential operators $\varphi (\Delta_{\theta , \omega})$, where $\varphi $
is an exponential type entire function of a single complex variable and $%
\Delta_{\theta , \omega} = (\theta + \omega z) D + zD^2$, $D = \partial /
\partial z$, $z\in \mathbb{C}$, $\theta \geq 0$, $\omega \in \mathbb{R}$,
acting in the spaces of exponential type entire function are studied. It is shown that,
for $\omega \geq0$, such operators preserve the set of Laguerre entire functions provided
the function $\varphi$ also belongs to this set. The latter consists of the polynomials
possessing real nonpositive zeros
only and of their uniform limits on compact subsets of the complex plane $%
\mathbb{C}$. The operator $\exp(a \Delta_{\theta , \omega})$, $a\geq 0$ is studied in
more details. In particular, it is shown that it preserves the set of Laguerre entire
functions for all $\omega \in \mathbb{R}$. An integral representation of $\exp(a
\Delta_{\theta ,\omega})$, $a>0$ is obtained. These results are used to obtain the
solutions to certain Cauchy problems employing $\Delta_{\theta , \omega}$.
\end{abstract}

\end{frontmatter}

\section{Introduction}

\noindent Differential operators of infinite order naturally appear in many applications
(in a certain sense they constitute a total set of linear operators acting between spaces
of differentiable functions \cite{Kor}). As usual, such operators are constructed by
means of finite order differential expressions substituted in the arguments of
appropriate functions. If such a function admits power series expansion, the operator may
be defined by imposing corresponding convergence conditions. In this article we consider
differentiation with respect to a single complex variable $z \in \mathbb{C}$ and the
operators constructed by means of the following differential expression
\begin{equation}  \label{1.1}
\Delta_{\theta , \omega} = \Delta_\theta + \omega z D \ \stackrel{\mathrm{def}%
}{=} \ (\theta + z D ) D + \omega z D,
\end{equation}
where $D = \partial/ \partial z$ and $\theta \geq 0$, $\omega \in \mathbb{R}$ are
parameters. Given entire functions $\varphi , f : \mathbb{C} \rightarrow \mathbb{C}$, we
set
\begin{equation}  \label{1.2}
\left( \varphi (\Delta_{\theta , \omega})f\right) (z) = \sum_{k=1}^\infty \frac{1}{k!}
\varphi^{(k)}(0) \left(\Delta_{\theta , \omega}^k f\right)(z).
\end{equation}
In order for the above series to converge to an entire function, we impose growth
restrictions on the functions $\varphi$ and $f$ by placing them into certain spaces of
exponential type entire function. These spaces were introduced and studied in our recent
work \cite{KozW1} where we used them to describe the operators $\varphi (\Delta_\theta)$,
$\theta \geq 0$. Our present research is mainly based on the results of that paper.

In Section 2 below we give definitions and a number of facts regarding the operators
$\varphi (\Delta_\theta)$, $\theta \geq 0$. Our main results are presented in Section 3
in a sequence of theorems characterizing the properties of $\varphi(\Delta_{\theta ,
\omega})$. Theorem \ref{2.1tm} describes the operator $\exp(a \Delta_{\theta , \omega})$.
In particular, we prove that this operator obeys the decomposition rule
\begin{equation}  \label{dec}
\exp(a\Delta_{\theta, \omega}) = \exp(a \omega z D)\cdot \exp\{ \omega^{-1} (e^{a \omega}
- 1)\Delta_{\theta}\},
\end{equation}
which is then used to describe the operator $\varphi (\Delta_{\theta,\omega})$ with an
arbitrary exponential type entire function $\varphi$ (Theorem \ref{2.2tm}). A special
role in our study is played by Laguerre entire functions being the polynomials of a
single complex variable possessing real nonpositive zeros only or the limits of sequences
of such polynomials taken in the topology of uniform convergence on compact subsets of
$\mathbb{C}$. It turns out that these functions are of exponential type and possess
corresponding infinite-product representations. We prove (Theorem \ref{2.3tm}) that if
both $\varphi$ and $f$ are Laguerre entire functions and if $\omega \geq 0$, the
$\varphi\left(\Delta_{\theta ,\omega}\right)f$ is also a Laguerre entire function. The
above theorems extend the results of \cite{KozW1} to nonzero values of $\omega$. Then we
consider again the operator $\exp\left(a\Delta_{\theta , \omega}\right)$, $%
a\geq 0$, for which the decomposition (\ref{dec}) implies that it preserves the set of
Laguerre entire functions \textit{for all real} $\omega$. Further, (\ref{dec}) is used to
obtain (Theorem \ref{2.5tm}) an integral representation of $\exp\left(a\Delta_{\theta ,
\omega}\right)$, $a > 0 $. The latter result allows us to extend this operator to a wider
class of functions. Theorem \ref{2.6tm} gives additional information regarding the action
of $\exp\left(a\Delta_{\theta, \omega} \right)$ on the functions of the type of $\exp(uz)
g (z)$. In Section 4 we use the above results to describe the solutions to the Cauchy
problem
\begin{eqnarray*}
\frac{\partial f (t , z)}{\partial t} & = & (\theta + \omega z) \frac{%
\partial f (t , z)}{\partial z} + z \frac{\partial^2 f (t , z)}{\partial z^2}%
, \quad t>0, \\
f (0, z ) & = & g(z),
\end{eqnarray*}
(Theorem \ref{e1tm}). It enables us to describe (Theorem \ref{Lunatm}%
) the solutions also to the Cauchy problem of the following type (diffusion equation with a
drift)
\begin{eqnarray*}
\frac{\partial F (t , x)}{\partial t} & = & \left(\Delta + b (x , \nabla)
\right) F(t , x), \quad t>0 , \\
F (0, x ) & = & G (x), \quad x \in \mathbb{R}^N, \ \ N \in \mathbb{N},
\end{eqnarray*}
with the initial function $G$ belonging to the class of isotropic (i.e. $O(N)
$-invariant) analytic functions. Here $\Delta$ and $\nabla$ stand for the $N$-dimensional
Laplacian and gradient respectively.

All Propositions are given without proofs since either they are taken from other sources
or the proofs are evident.

\section{Preliminaria}

Let $\mathcal{E}$ stand for the set of all entire functions $f : \mathbb{C} \rightarrow
\mathbb{C}$ equipped with the topology $\mathcal{T}_{\mathbb{C}}$ of uniform convergence
on
compact subsets of $\mathbb{C}$. Thus, $(\mathcal{E}, \mathcal{T}_{\mathbb{C}%
})$ is a Fr{\'e}chet space. For $b>0$, we set
\[
\mathcal{B}_{b}=\{f\in {\mathcal{E}}\mid \Vert f\Vert _{b}<\infty \};
\]
where
\begin{equation}  \label{1}
\Vert f\Vert _{b}=\sup_{k\in {\mathbb{N}}_{0}}\{b^{-k}\mid f^{(k)}(0)\mid \}, \ \ \
f^{(k)}(0)=(D^{k}f)(0),
\end{equation}
and $\mathbb{N}_{0}$ stands for the set of all nonnegative integers. Every $%
\left( \mathcal{B}_{b},\Vert \cdot \Vert _{b}\right) $ is a Banach space. For $a\geq 0$,
let
\begin{equation}
\mathcal{A}_{a}=\bigcap_{b>a}\mathcal{B}_{b}=\{f\in \mathcal{E}\mid (\forall b>a)\ \Vert
f\Vert _{b}<\infty \}.  \label{5}
\end{equation}
Equipped with the topology $\mathcal{T}_{a}$ defined by the family of norms $\{\Vert
.\Vert _{b},\ b>a\}$, this set becomes a Fr\'{e}chet space. To shorten our notation we
write $\mathcal{E}$, $\mathcal{A}_a $, $\mathcal{B}_b $ instead of $(\mathcal{E},
\mathcal{T}_{\mathbb{C}})$, $(\mathcal{A}_a , \mathcal{T}_a)$, $(\mathcal{B}_b , \| \cdot
\|_b )$ respectively.

\begin{Df}
\label{1df} A family $\mathcal{L}$ (respectively $\mathcal{L}_{0}$, $%
\mathcal{\ L}^{+}$, $\mathcal{L}^{-}$) consists of the entire functions possessing the
representation
\begin{eqnarray}  \label{6}
& & f(z)=Cz^{m}\exp (\alpha z)\prod_{j=1}^{\infty }(1+\beta _{j}z); \\
& & C\in\mathbb{C} ,\ m\in \mathbb{N}_{0},\ \ \beta _{j}\geq \beta _{j+1}\geq 0, \
\sum_{j=1}^{\infty }\beta _{j}<\infty ,  \nonumber
\end{eqnarray}
with $\alpha \in \mathbb{R}$ (respectively $\alpha =0$, $\alpha \geq 0$, and $\alpha
<0$).
\end{Df}
The elements of ${\mathcal{L}^{+}}$ are known as Laguerre entire functions \cite{Il},
\cite{KozW1}. Let $\mathcal{P}^{+}$ stand for the set of polynomials belonging to
$\mathcal{L}^+$. By Laguerre and P{\'o}lya (see e.g. \cite{Il}, \cite{Lev}),
\begin{Pn}
\label{4pn} The family $\mathcal{L}^{+}$ is exactly the closure of $\mathcal{ P}^+$ in
$\mathcal{T}_{\mathbb{C}}$.
\end{Pn}
It is worth to note that every $f$ being of the form (\ref{6}) may be written $f(z)=\exp
(\alpha z)h(z)$, where $h$ is an entire function of order less than one or equal to one
and, in the latter case, of minimal type. Consider the families
\begin{equation}
\mathcal{L}_{a}\ {\stackrel{\mathrm{def}}{=}}\ \mathcal{L}\cap \mathcal{A}%
_{a},\ \ \mathcal{L}_{a}^{\pm }\ {\stackrel{\mathrm{def}}{=}}\ \mathcal{L}%
^{\pm }\cap \mathcal{A}_{a}.  \label{7}
\end{equation}
Obviously, $\mathcal{P}$ and $\mathcal{P}^{+}$ are dense respectively in the sets
$\mathcal{A}_{a}$ and $\mathcal{L}_{a}^{+}$ equipped with the topologies induced on them
by $\mathcal{T}_\mathbb{C}$. However, the set $\mathcal{P}$ is not dense in any space
$\mathcal{B}_{b}$.
Thus, \textit{a priori} it is not obvious whether or not the sets $\mathcal{P%
}$ and $\mathcal{P}^{+}$ are dense respectively in $\mathcal{A}_{a}$ (in its
standard topology) and in $\mathcal{L}_{a}^+$ in the topology induced from $%
\mathcal{A}_{a}$. Fortunately, this density property holds in both cases. The following
two statements, which we borrow from \cite{KozW1}, give
information regarding the topological properties of $\mathcal{A}_a$ and $%
\mathcal{L}^+_a$.

\begin{Pn}
\label{T00} For every $a\geq 0,$ the relative topology on a bounded subset
of $\mathcal{A}_{a}$ coincides with the topology induced on it by $\mathcal{T%
}_{\mathbb{C}}$.
\end{Pn}

\begin{Pn}
\label{T0} For every $a\geq 0$,\\[.1cm]
\begin{tabular}{ll}
(i) & the set of all polynomials $\mathcal{P}\subset
\mathcal{E}$ is dense in $\mathcal{A}_{a}$; \\[0.1cm]  (ii) &
the set $\mathcal{P}^{+}$ is dense in $\mathcal{L}_{a}^{+}$ in the topology induced by
$\mathcal{T}_{a}$.
\end{tabular}
\end{Pn}
Below unless explicitly stated we consider $\mathcal{L}^+_a$, $a\geq 0$ as a topological
space equipped with the topology induced by $\mathcal{T}_a$.

For $\theta \geq 0$, let $\Delta _{\theta }:\mathcal{E} \rightarrow \mathcal{E}$ be as in
(\ref{1.1}), that is $\Delta _{\theta } = \Delta_{\theta , 0} = (\theta + z D)D$. Given
entire functions $\varphi $ and $f$, we define $\varphi (\Delta _{\theta })f(z)$ by
(\ref{1.2}). Further, one has
\begin{eqnarray}  \label{q}
\Delta _{\theta }^{k}z^{m} & = & q_{\theta }^{(m,k)}z^{m-k}, \\
q_{\theta }^{(m,k)} & =& \left\{
\begin{array}{ll}
0, & k>m \\
&  \\
\gamma _{\theta }(m)/\gamma _{\theta }(m-k), & 0\leq k\leq m
\end{array}
\right. ,  \nonumber
\end{eqnarray}
where
\[
\gamma _{\theta }(m)=m!\Gamma (\theta +m).
\]
Applying this in (\ref{1.2}) one may prove the following statement \cite{KozW1}.

\begin{Pn}
\label{PrB} For all $\theta \geq 0$ and for arbitrary $a>0$ and $b>0$, such that $ab<1$,
$(\varphi ,f)\mapsto \varphi (\Delta _{\theta })f$ is a continuous bilinear map from
$\mathcal{B}_{a}\times \mathcal{B}_{b}$ (resp. from $\mathcal{A} _{a}\times
\mathcal{A}_{b}$ with $a\geq 0,$ $b\geq 0$) into $\mathcal{B}_{c}$ (resp.
$\mathcal{A}_{c}$), where $c=b(1-ab)^{-1}$. Moreover,
\[
\|\varphi(\Delta_{\theta})f\|_c \leq (1- ab)^{-\theta}\|\varphi\|_a \|f\|_b .
\]
\end{Pn}

The action of $\varphi (\Delta_\theta )$ on the Laguerre entire functions is described by
the following statement, which was proven in \cite{KozW1}.

\begin{Pn}
\label{a1tm} For all $\theta \geq 0$ and for arbitrary $a\geq 0$ and $b\geq 0 $, such
that $ab<1$, $(\varphi ,f)\mapsto \varphi (\Delta _{\theta })f$ is a
continuous map from $\mathcal{L}_{a}^{+}\times \mathcal{L}_{b}^{+}$ into $%
\mathcal{L} _{c}^{+},$ where $c=b(1-ab)^{-1}$.
\end{Pn}

Given $\omega \in \mathbb{R}$ and $f\in \mathcal{E}$, we set
\[
\exp (\omega z D)f(z)=\sum_{k=0}^{\infty }\frac{\omega^{k}}{k!}%
\left((zD)^{k}f\right)(z),
\]
which readily yields
\begin{equation}  \label{3.1}
\exp (\omega zD)f(z)=f(e^{\omega}z).
\end{equation}

\begin{Pn}
\label{pro5} For any $\omega \in \mathbb{R}$, $\exp (\omega zD)$ is a
continuous linear map from $\mathcal{B}_{c}$ with $c>0$ (resp. from $%
\mathcal{A}_{c}$ with $c\geq 0$) into $\mathcal{B}_{d}$ (resp. $\mathcal{A}%
_{d}$), where $d=e^{\omega}c$.
\end{Pn}

The identity (\ref{3.1}) also implies that

\begin{Pn}
\label{Pro6} For any $\omega \in \mathbb{R}$ and $c\geq 0,$ $\exp (\omega zD)
$ is a continuous map from $\mathcal{L}_{c}^{+}$ (resp. from $\mathcal{L}%
_{c}^{-}$) into $\mathcal{L}_{d}^{+}$ (resp. $\mathcal{L}_{d}^{-}$ ), where $%
d=e^{\omega}c$.
\end{Pn}

The following representation of $\exp(a \Delta_\theta )$ was obtained in \cite{KozW1}.

\begin{Pn}
\label{intt1pn} For every $\theta \geq 0$, for arbitrary $a>0$, $b \geq 0$, such that
$ab<1$, and for all $f \in \mathcal{A}_b$,
\begin{eqnarray}  \label{intt1}
\left(\exp( a\Delta_\theta )f \right)(z) & = & \exp\left( -\frac{z}{a} \right)
\int_{0}^{+\infty} w_\theta \left(\frac{sz}{a}\right) f(as)
s^{\theta -1} e^{-s} \mathrm{d} s  \nonumber \\
& \stackrel{\mathrm{def}}{=} & \int_{0}^{+\infty} K_\theta \left(\frac{z}{a}, s\right)
f(as) s^{\theta -1} e^{-s} \mathrm{d} s ,
\end{eqnarray}
where
\begin{equation}  \label{int3}
K_\theta (z , s) = e^{-z} w_{\theta} (zs), \quad w_\theta (\xi ) \ \stackrel{%
\mathrm{def}}{=} \ \sum_{k=0}^\infty \frac{\xi^k }{\gamma_\theta (k)}.
\end{equation}
\end{Pn}

This representation may be used for extending the operator $\exp (a \Delta_\theta)$. The
assertion below, also taken from \cite{KozW1}, describes a property of such an extended
operator.

\begin{Pn}
\label{intt2pn} Given $a>0$ and $u \in \mathbb{R}$, let $b \in [0, - u + 1/a )$. Then,
for every $g \in \mathcal{A}_b $, the operator (\ref{intt1}) may be applied to the
function
\begin{equation}  \label{intt2}
f(z) = \exp ( u z) g(z),
\end{equation}
yielding
\begin{equation}  \label{intt3}
\left(\exp( a \Delta_\theta ) f \right) (z) = (1 - ua)^{-\theta} \exp\left(%
\frac{u z }{1 - ua} \right) h(z),
\end{equation}
where
\begin{eqnarray}  \label{intt4}
h(z) & = & \left[\exp\left( \frac{a}{1-ua} \Delta_\theta\right)g \right]
\left( \frac{z}{(1- ua)^2}\right)  \nonumber \\
& = & \exp\left[a(1-ua) \Delta_\theta \right]g \left(\frac{z}{(1- ua)^2} \right).
\end{eqnarray}
Moreover, $h \in \mathcal{A}_c$ with $c = b (1-ua)^{-1} [1 - a(u+b)]^{-1}$.
\end{Pn}

The extensions of $\exp(a\Delta_{\theta})$ on the base of the integral representation
(\ref{intt1}) to the spaces of integrable functions are given in \cite{KOU}.

\section{Main Results}

\begin{Tm}
\label{2.1tm} For all $\theta \geq 0$, for arbitrary $\omega \in \mathbb{R}$, $a\geq0$,
and $b >0$, obeying the condition $b \omega^{-1}[e^{a\omega} - 1] < 1$, the series
expansion (\ref{1.2}) defines a continuous linear operator $\exp\left( a\Delta_{\theta ,
\omega} \right)$
acting from the space $\mathcal{B}_b $ (resp. from $\mathcal{A}_b$ with $%
b\geq 0$) into the space $\mathcal{B}_c$ (resp. $\mathcal{A}_c$), where
\begin{equation}  \label{2.1}
c = c(b) \ \stackrel{\mathrm{def}}{=} \ b e^{a\omega}\left(1 - b \omega^{-1}[e^{a \omega}
- 1]\right)^{-1}.
\end{equation}
This operator obeys the decomposition rule (\ref{dec}).
\end{Tm}

\begin{Rk}
\label{R1} This theorem holds also for $\omega=0,$ since $\lim_{\omega \rightarrow
0}[e^{a \omega }-1]\omega^{-1} =a,$ where it coincides with a partial case of Proposition
\ref{PrB}. In the sequel the case $\omega=0$ will be understood in this sense.
\end{Rk}

\begin{Tm}
\label{2.2tm} For all $\theta \geq 0$, for arbitrary $a > 0$, $\omega \in \mathbb{R}$, $b
>0$, obeying the condition $b \omega^{-1}[e^{a \omega} - 1] < 1$, $(\varphi , f) \mapsto
\varphi(\Delta_{\theta , \omega})f $, is a continuous bilinear map from $\mathcal{B}_a
\times \mathcal{B}_b $ (resp. from $\mathcal{A}_a \times \mathcal{A}_b$ with $a \geq 0$
and $b \geq 0$) into the space $\mathcal{B}_c$ (resp. $\mathcal{A}_c$), where $c = c(b)$
is given by (\ref{2.1}). Moreover,
\[
\|\varphi(\Delta_{\theta , \omega})f\|_c \leq \left(1 - b \omega^{-1}[e^{a \omega} -1]
\right)^{-\theta} \|\varphi\|_a \|f\|_b .
\]
\end{Tm}

\begin{Tm}
\label{2.3tm} For all $\theta \geq 0$, for arbitrary $\omega \geq 0$, $a \geq 0$, $b \geq
0$, obeying the condition $b \omega^{-1}[e^{a\omega} - 1] < 1$, $(\varphi, f) \mapsto
\varphi(\Delta_{\theta , \omega}) f$ is a
continuous map from $\mathcal{L}^+_a \times \mathcal{L}^+_b$ into $\mathcal{L%
}^+_c$, where $c = c(b)$ is given by (\ref{2.1}).
\end{Tm}

\begin{Tm}
\label{2.30tm} For all $\theta \geq 0$, for arbitrary $\omega \in \mathbb{R}$%
, $a \geq 0$, $b \geq 0$, obeying the condition $b \omega^{-1}[e^{a \omega}
- 1] < 1$, $\exp(a\Delta_{\theta, \omega}) $ is a continuous map from $%
\mathcal{L}^+_b$ into $\mathcal{L}^+_c$, where $c = c(b)$ is given by (\ref {2.1}).
\end{Tm}

\begin{Tm}
\label{2.5tm} For all $\theta \geq 0$, for arbitrary $\omega \in \mathbb{R}$%
, $a > 0$, and $b\geq 0$, obeying the condition $b \omega^{-1}[e^{a \omega} - 1] < 1$,
and for any $f \in \mathcal{A}_b$,
\begin{eqnarray}  \label{int1}
& & \left(\exp(a \Delta_{\theta , \omega})f \right) (z) = \\
& & \quad = \exp\left(- \nu z \right) \int_{0}^{+\infty}w_{\theta} (\nu s z) f
(\omega^{-1}[e^{a \omega } -1]s) s^{\theta -1} e^{-s} \mathrm{d} s  \nonumber
\\
& & \quad = \int_{0}^{+\infty} K_{\theta} (\nu z, s) f (\omega^{-1}[e^{a \omega } -1]s)
s^{\theta -1} e^{-s} \mathrm{d} s ,  \nonumber
\end{eqnarray}
where
\begin{equation}  \label{int2}
\nu = \nu( a , \omega ) \ \stackrel{\mathrm{def}}{=} \ \frac{\omega e^{a \omega}}{e^{a
\omega } - 1},
\end{equation}
and $w_\theta $ and $K_\theta$ are defined by (\ref{int3}).
\end{Tm}

The above theorem allows us to extend the operator $\exp( a \Delta_{\theta , \omega})$ to
the functions for which the integrals in the right-hand side of (\ref{int1}) converge. In
the sequel we understand this operator in such an extended (integral) version.

\begin{Tm}
\label{2.6tm} Given $a \geq 0$, $\omega \in \mathbb{R}$, $u \in \mathbb{R}$, and $b \geq
0$, obeying the condition $b \omega^{-1}[e^{a \omega} - 1] < 1 - u \omega^{-1}[e^{a
\omega} - 1] $, the operator $\exp(a \Delta_{\theta , \omega})$ may be applied to the
function
\begin{equation}  \label{int4}
f(z) = \exp( u z) g(z),
\end{equation}
with arbitrary $g \in \mathcal{A}_b$, yielding
\begin{eqnarray}  \label{int5}
& & \left(\exp(a \Delta_{\theta , \omega})f \right) (z) = \\
& & \quad = \left(1 - u \omega^{-1}[e^{a \omega} - 1] \right)^{-\theta} \exp\left(\frac{u
e^{a \omega}z}{1 - u \omega^{-1}[e^{a \omega} - 1]} \right) h(z),  \nonumber
\end{eqnarray}
where
\begin{eqnarray}  \label{int6}
h(z) & = & \left[ \exp\left(\frac{a\omega}{\omega - u[e^{a \omega} - 1]}
\Delta_\theta + a \omega z D\right)g \right]\left(\frac{\omega^2 z} {%
\left(\omega - u [e^{a \omega} - 1]\right)^2} \right)  \nonumber \\
& = & \left[ \exp\left(\frac{e^{a \omega} - 1}{\omega - u [e^{a \omega} - 1]}
\Delta_\theta \right)g \right]\left( \frac{\omega^2 e^{a \omega}z} {%
\left(\omega - u [e^{a \omega} - 1]\right)^2} \right).
\end{eqnarray}
Moreover, the latter function belongs to $\mathcal{A}_c$ with
\[
c = b e^{a \omega} \left( 1 - u \omega^{-1}[e^{a \omega} - 1]\right)^{-1} \left( 1-
(u+b)\omega^{-1}[e^{a \omega} - 1]\right)^{-1}.
\]
\end{Tm}

\begin{Co}
\label{2.1co} For all $\theta \geq 0$, $\omega \in \mathbb{R}$, $a\geq 0$, and $b\geq 0$,
the operator $\exp\left( a \Delta_{\theta , \omega}\right)$ is a continuous map from
\newline
\begin{tabular}{ll}
\textrm{(i)} & $\mathcal{L}_b^+ $ into $\mathcal{L}_c^+$, where $a$, $b$,
and $\omega$ satisfy the condition \\
& $b\omega^{-1}[e^{a\omega} - 1] < 1$, and $c= c(b)$ is given by (\ref{2.1}),
\\[.1cm]
\textrm{(ii)} & $\mathcal{L}_b^{-}$ into $\mathcal{L}_d^{-}$, where $d = b e^{a \omega}
\left( 1+b \omega^{-1}[e^{a \omega} - 1]\right)^{-1}$.
\end{tabular}
\end{Co}

The proof of the above theorems will be based on the properties of $\exp(a \Delta_{\theta
, \omega})$, which may be studied on the base of (\ref{dec})
and the properties of $\exp(a \omega z D)$ and $%
\exp(b \omega^{-1} [e^{a \omega} -1] \Delta_\theta )$, given by Propositions \ref{PrB} -
\ref{intt2pn}. Clearly, the operator $\exp(a\Delta_{\theta ,\omega})$ defined by
(\ref{1.2}) may be applied to any $f \in \mathcal{P}$, furthermore, $\exp(a
\Delta_{\theta ,\omega}): \mathcal{P} \rightarrow \mathcal{P}$. Similarly, for any real
$b$ and $c$, $\exp(b z D): \mathcal{P} \rightarrow \mathcal{P}$ and $\exp(c \Delta_\theta
): \mathcal{P} \rightarrow \mathcal{P}$.

\begin{Lm}
\label{main} For every polynomial $f\in \mathcal{P}$ and any $a\in \mathbb{R} $, the
operator $\exp(a \Delta_{\theta , \omega})$ defined by the series expansion (\ref{1.2})
obeys (\ref{dec}).
\end{Lm}

\textbf{Proof.} Obviously, it is enough to prove the lemma only for $f_m (z) = z^m$,
$m\in \mathbb{N}$. Since $a \Delta_\theta $ and $a \omega z D$ do not commute (except for
$a$ or $\omega$ being zero), one has to apply an appropriate decomposition technique. It
will be based on the Trotter-Kato product formula, which, for our aims, may be written as
follows (for more details on this item see \cite{Nzag} and the references therein)
\begin{equation}  \label{TK}
\exp (A+B)f=\lim_{n\rightarrow \infty }\left( \exp \left( \frac{1}{n}%
A\right) \exp \left( \frac{1}{n}B\right) \right) ^{n}f,
\end{equation}
where $A$ and $B$ are linear continuous operators on a topological vector space and the
convergence is understood in the topology of the range space.

In what follows, given $m\in \mathbb{N}$, one has
\begin{equation}
\exp (a\Delta _{\theta ,\omega })f_{m}(z)=\lim_{n\rightarrow \infty }\left(
\exp \left( \frac{a}{n}\Delta _{\theta }\right) \exp \left( \frac{a\omega }{n%
}zD\right) \right) ^{n}f_{m}(z),  \label{tr1}
\end{equation}
where the convergence is in $\mathcal{T}_{\mathbb{C}} $. By (\ref{pro5}) and (\ref{q}),
\[
\exp (a{\Delta }_{\theta })z^{m}=\sum_{k=0}^{m}\frac{a^{k}}{k!}q_{\theta
}^{(m,k)}z^{m-k},
\]
which yields, together with (\ref{3.1}), that for any $n,m\in \mathbb{N}$,
\begin{eqnarray}
&&\qquad \left( \exp \left( \frac{a}{n}\Delta _{\theta }\right) \exp \left(
\frac{a\omega }{n}zD\right) \right) ^{n}z^{m}=  \nonumber  \label{a} \\
&=&\sum_{k_{1}=0}^{m}\sum_{k_{2}=0}^{m-k_{1}}\dots
\sum_{k_{n}=0}^{m-k_{1}-\dots -k_{n-1}}\frac{1}{k_{1}!k_{2}!\dots k_{n}!}%
\times   \nonumber \\
&\times &\left( \frac{a}{n}\right) ^{k_{1}+k_{2}+\dots +k_{n}}q_{\theta
}^{(m,k_{1})}q_{\theta }^{(m-k_{1},k_{2})}\dots q_{\theta }^{(m-k_{1}-\dots
-k_{n-1},k_{n})}\times   \nonumber \\
&\times &\exp \left\{ \frac{a\omega }{n}(m+m-k_{1}+\dots +m-k_{1}-\dots
-k_{n-1})\right\} z^{m-k_{1}-\dots -k_{n}}  \nonumber \\
&=&\sum_{k_{1}=0}^{m}\sum_{k_{2}=0}^{m-k_{1}}\dots
\sum_{k_{n}=0}^{m-k_{1}-\dots -k_{n-1}}\frac{1}{k_{1}!k_{2}!\dots k_{n}!}%
\times   \nonumber \\
&&\qquad \times \left( \frac{a}{n}\right) ^{k_{1}+k_{2}+\dots
+k_{n}}z^{m-k_{1}-\dots -k_{n}}\times   \nonumber \\
&&\qquad \times \frac{\gamma _{\theta }(m)\exp \left\{ \frac{a\omega }{n}%
(nm-(n-1)k_{1}-\dots -k_{n-1})\right\} }{\gamma _{\theta }(m-k_{1}-\dots
-k_{n})}  \nonumber \\
&&\qquad \stackrel{\mathrm{def}}{=}\ \sum_{k=0}^{m}a^{k}\frac{\gamma _{\theta
}(m)}{\gamma _{\theta }(m-k)}\left( e^{a\omega }z\right) ^{m-k}\Psi _{n}(k,a\omega ).
\end{eqnarray}
Here for $k\in \mathbb{N}_{0}$ and $b\in \mathbb{R}$, we have set
\begin{eqnarray}
\Psi _{n}(k,b) &=&\frac{1}{n^{k}}\sum_{k_{1},k_{2},\dots ,k_{n}=0}^{k}\delta
_{k,k_{1}+\dots +k_{n}}\frac{e^{\frac{b}{n}(k_{1}+2k_{2}+\dots +nk_{n})}}{%
k_{1}!k_{2}!\dots k_{n}!}  \nonumber  \label{b} \\
&=&\frac{1}{k!n^{k}}\left( e^{\frac{b}{n}}+e^{2\frac{b}{n}}+\dots +e^{n\frac{%
b}{n}}\right) ^{k}=\frac{1}{k!}\left( \frac{e^{\frac{b}{n}}\left(
e^{b}-1\right) }{n\left( e^{\frac{b}{n}}-1\right) }\right) ^{k}  \nonumber \\
&\longrightarrow &\frac{1}{k!}\left( \frac{e^{b}-1}{b}\right) ^{k},
\end{eqnarray}
when $n\rightarrow +\infty $. By (\ref{tr1}), (\ref{a}), and (\ref{b}) one obtains
\begin{eqnarray*}
&&\exp (a\Delta _{\theta ,\omega })f_{m}(z)=\lim_{n\rightarrow \infty
}\left( \exp \left( \frac{a}{n}\Delta _{\theta }\right) \exp \left( \frac{%
a\omega }{n}zD\right) \right) ^{n}z^{m} \\
&&\quad =\sum_{k=0}^{m}\frac{\left( \omega ^{-1}[e^{a\omega }-1]\right) ^{k}%
}{k!}\cdot \frac{\gamma _{\theta }(m)}{\gamma _{\theta }(m-k)}\left(
e^{a\omega }z\right) ^{m-k} \\
&&\quad =\left( \exp \left( \omega ^{-1}[e^{a\omega }-1]{\Delta }_{\theta
}\right) f_{m}\right) (e^{a\omega }z) \\
&&\quad =\exp (a\omega zD)\exp (\omega ^{-1}[e^{a\omega }-1]\Delta _{\theta })f_{m}(z).
\end{eqnarray*}
\quad \hfil%
\hbox{
    \vrule width 0.8pt height 7pt depth 0.8pt  \hskip -0.2pt
    \vrule width   5pt height 7pt depth -6.2pt \hskip -5pt
    \vrule width   5pt height 0pt depth 0.8pt
    \vrule width 0.8pt height 7pt depth 0.8pt \hskip 2pt}

Given $k,n,m\in \mathbb{N}$, we set
\begin{equation}
\kappa _{n}(k,m,\theta ,\omega )=\left[ D^{n}\left( \Delta _{\theta ,\omega }\right)
^{k}z^{m}\right] _{z=0}.  \label{kappa}
\end{equation}
After some algebra one obtains
\begin{eqnarray}
&&\kappa _{n}(k,m,\theta ,\omega )=  \label{kappa1} \\
&&\qquad =\frac{n!}{(m-n)!}q_{\theta }^{(m,m-n)}\omega ^{k+n-m}\sum_{l=m-n}^{k}{{k
\choose l}}n^{k-l}\alpha _{m-n}^{(l)},  \nonumber
\end{eqnarray}
if $m\geq n$, and $\kappa _{n}(k,m,\theta ,\omega )=0$ if $m<n$. Here, for $%
p\geq s$, $p,s\in \mathbb{N}$,
\[
\alpha _{p}^{(s)}\stackrel{\mathrm{def}}{=}\sum_{k=0}^{p}{{p \choose k}}%
(-1)^{k}(p-k)^{s}>0.
\]

\textbf{Proof of Theorem \ref{2.1tm}.} Given $f \in \mathcal{P}$ and $b >0$
obeying the condition $b\omega^{-1} [e^{a \omega} -1] < 1$, one has from (%
\ref{dec}), (\ref{3.1}), and Proposition \ref{PrB}
\begin{eqnarray}  \label{est}
\|\exp( a\Delta_{\theta , \omega}) f \|_c & = & \|\exp(a \omega z D ) \exp\{
\omega^{-1}(e^{a \omega}-1)\} f\|_c =  \nonumber \\
& = & \|\exp\{ \omega^{-1}(e^{a \omega}-1)\} f\|_{ce^{-a\omega}}  \nonumber \\
&\leq & \left( 1 - b \omega^{-1} (e^{a\omega } - 1)\right)^{-\theta} \|f\|_b ,
\end{eqnarray}
where $c = c(b)$ is given by (\ref{2.1}). Then by claim (i) of Proposition
\ref{T0},$\exp(a \Delta_{\theta , \omega})$ may be continuously extended to the whole
$\mathcal{A}_b$. The extension will also obey (\ref{dec}). This yields in turn that the
estimate (\ref{est}) holds for any $f \in \mathcal{E}$, provided $\|f\|_b <
\infty$.  \quad\hfil%
\hbox{
    \vrule width 0.8pt height 7pt depth 0.8pt  \hskip -0.2pt
    \vrule width   5pt height 7pt depth -6.2pt \hskip -5pt
    \vrule width   5pt height 0pt depth 0.8pt
    \vrule width 0.8pt height 7pt depth 0.8pt \hskip 2pt}

\textbf{Proof of Theorem \ref{2.2tm}.} According to (\ref {1.2}), (\ref{kappa})
\begin{eqnarray}  \label{gamma2}
g^{(n)}(0) & \stackrel{\mathrm{def}}{=} & \left(\varphi (\Delta_{\theta ,
\omega})f\right)^{(n)} (0) \\
& = & \sum_{k, m =0}^{\infty} \frac{\varphi^{(k)}(0)}{k!} \cdot \frac{%
f^{(m)}(0)}{m!} \kappa_n (k , m, \theta , \omega), \quad n \in \mathbb{N}_0 , \nonumber
\end{eqnarray}
which, for given positive $a $ and $b$, yields
\begin{equation}  \label{gamma22}
|g^{(n)} (0)| \leq \|\varphi \|_a \|f\|_b \sum_{k, m =0}^{\infty} \frac{a^k}{%
k!} \cdot \frac{b^m}{m!}| \kappa_n (k , m, \theta , \omega)|.
\end{equation}
By (\ref{kappa1}),
\[
| \kappa_n (k , m, \theta , \omega)| \leq \kappa_n (k , m, \theta , |\omega|).
\]
Therefore, for $a$ and $b$ obeying the conditions of the theorem and for $c (b)$ given by
(\ref{2.1}), one readily gets from (\ref{gamma22}), (\ref {gamma2}), (\ref{est}), and
(\ref{1})
\begin{eqnarray*}
\|g \|_c & \leq & \|\varphi \|_a \|f \|_b \cdot \sup_{n \in \mathbb{N}_0} \left\{c^{-n}
\sum_{k, m =0}^{\infty} \frac{a^{k}}{k!} \cdot \frac{b^{m}}{m!}
\kappa_n (k , m, \theta , |\omega|)\right\} \\
& = & \|\varphi \|_a \|f \|_b \cdot \| \exp\{ a \Delta_{\theta , |\omega|}\}
\exp(b z) \|_c \\
& \leq & \left( 1 - b\frac{e^{a |\omega|} - 1}{|\omega|}\right)^{-\theta}\|%
\varphi\|_a\|f\|_b = \left( 1 - b\frac{e^{a \omega} - 1}{\omega}%
\right)^{-\theta}\|\varphi\|_a\|f\|_b.
\end{eqnarray*}
\quad\hfil%
\hbox{
    \vrule width 0.8pt height 7pt depth 0.8pt  \hskip -0.2pt
    \vrule width   5pt height 7pt depth -6.2pt \hskip -5pt
    \vrule width   5pt height 0pt depth 0.8pt
    \vrule width 0.8pt height 7pt depth 0.8pt \hskip 2pt}

\begin{Lm}
\label{l1lm} For all $\theta \geq 0$, $\omega \geq 0$ and for arbitrary $%
\varphi, f \in \mathcal{P}^+$, the polynomial $\varphi (\Delta_{\theta , \omega}) f$ also
belongs to $\mathcal{P}^+$.
\end{Lm}

This lemma will be proven in several steps below. The proof of Theorem \ref {2.3tm}
readily follows from it and from claim (ii) of Proposition \ref{T0} and Theorem
\ref{2.2tm}. Set
\begin{equation}  \label{l1}
U = \{ z \in \mathbb{C} \ | \ \mathrm{Re}z > 0\}.
\end{equation}
In \cite{KozW1} we have proven the following statement.

\begin{Pn}
\label{l2pn} Let $P$, $Q$, and $Q_1$ be polynomials of a single complex variable. Suppose
that $P$ does not vanish on $U$ and
\begin{equation}  \label{l2}
Q (u) + P(u) v Q_1 (u) \neq 0 ,
\end{equation}
whenever $u, v \in U$. Then either
\begin{equation}  \label{l3}
S(z) \ \stackrel{\mathrm{def}}{=} \ Q(z) + P(z) D Q_1(z) \neq 0,
\end{equation}
whenever $z \in U$, or else $S(z) \equiv 0$.
\end{Pn}

\begin{Lm}
\label{l2lm} For arbitrary $\sigma \geq 0$, $\theta \geq 0$, and $\omega \geq 0$, the
second order differential operator $\sigma + \Delta_{\theta , \omega}$ maps
$\mathcal{P}^+$ into itself.
\end{Lm}

\textbf{Proof.} For an arbitrarily chosen $p\in \mathcal{P}^+$, we have to show that
$(\sigma + \Delta_{\theta , \omega})p\in \mathcal{P}^+$. The case of constant $p$ is
trivial. For nonconstant $p\in \mathcal{P}^+$, one may write
\begin{equation}  \label{l4}
p(z) = \pi_0 \prod_{j=1}^m (\pi_j + z), \ \ m\in \mathbb{N}, \ \ \pi_j \geq 0, \ \ j = 1,
2, \dots , m.
\end{equation}
First we consider the case $\pi_j = 0$, $j = 1 , 2 , \dots , m$, that is $%
p(z) = \pi_0 z^m$. Then
\[
(\sigma + \Delta_{\theta , \omega})p(z) = z^{m-1} \tilde{p}(z) \in \mathcal{P%
}^+,
\]
since
\[
\tilde{p}(z) \ \stackrel{\mathrm{def}}{=}m \ (\sigma + m \omega )z + m ( \theta + m-1)
\in \mathcal{P}^+.
\]
Now let at least one $\pi_j $ in (\ref{l4}) do not vanish. Set
\begin{equation}  \label{l6}
q(z) = p(z^2 ) = \pi_0 \prod_{j=1}^m (\pi_j + z^2 ).
\end{equation}
Then
\begin{equation}  \label{l7}
\left( \left( \sigma + \Delta_{\theta , \omega}\right)p \right)(z^2) = \left( \left(
\sigma + \Lambda_{\theta , \omega}\right)q \right)(z),
\end{equation}
where
\begin{equation}  \label{l8}
\Lambda_{\theta , \omega} \ \stackrel{\mathrm{def}}{=} \ \left( \theta +
\frac{z}{2}D\right)\left(\frac{1}{2z} D \right) + \frac{\omega}{2} z D.
\end{equation}
In view of (\ref{l6}), the polynomials $q(z)$, $q(-z)$ do not vanish on $U$ (\ref{l1}).
The proof will be done by showing that
\begin{equation}  \label{l80}
\left( \left( \sigma + \Lambda_{\theta , \omega}\right)q \right)(z) \neq 0,
\end{equation}
whenever $z \in U$. Taking into account (\ref{l6}) we may write
\[
\left( \left( \sigma + \Lambda_{\theta , \omega}\right)q \right)(z) = [ \sigma + ( \theta
+ \omega z^2 )r(z) ] q(z) + \frac{1}{2} z D \left(q(z) r(z) \right),
\]
where
\[
r(z) = \sum_{j=1}^{m} \frac{1}{\pi_j + z^2}.
\]
Set
\begin{eqnarray*}
Q(z) & = & [ \sigma + ( \theta + \omega z^2 )r(z) ] q(z), \quad P(z) = \frac{%
z}{2}, \\
& & \quad Q_1 (z) = q(z) r(z) .
\end{eqnarray*}
By Proposition \ref{l2pn}, the proof of (\ref{l80}) will be done if we show that
\begin{equation}  \label{l12}
R( u , v) \ \stackrel{\mathrm{def}}{=} \ Q(u) + P(u) v Q_1 (u) \neq 0,
\end{equation}
whenever $u, v \in U$. To this end we rewrite the latter as follows
\begin{equation}  \label{l13}
R( u , v)= \frac{1}{2} q(u) \cdot R_1 (u , v) \cdot R_2 (u , v),
\end{equation}
where
\[
R_1 ( u , v) = 2 \omega u^2 + u v + 2 \theta ,
\]
and
\begin{equation}  \label{l15}
R_2 ( u , v) = r(u) + \frac{2 \sigma}{2 \omega u^2 + u v + 2 \theta}.
\end{equation}
Since $q(u) \neq 0$ whenever $u \in U$, to prove (\ref{l12}) it remains to
show that both $R_1$ and $R_2$ do not vanish if $u , v\in U$. Given $u\in U$%
, let us solve the equation $R_1 (u , v) = 0$. The result is
\[
v = - \frac{2\theta}{|u|^2} \bar{u} - 2 \omega u,
\]
which yields the following implications
\[
\left(R_1 (u , v) = 0 \right) \Rightarrow \left(\mathrm{Re} v \leq 0 \right) \Rightarrow
\left( v \in \mathbb{C} \setminus U \right).
\]
Recall that both $\theta $ and $\omega $ are supposed to be real and nonnegative. Then
$R_1 (u , v) \neq 0$ if $v \in
U$. By the same arguments we show that $R_2 (u , v) \neq 0$ if $u , v \in U$%
. To this end we rewrite (\ref{l15})
\begin{equation}  \label{l16}
R_2 (u , v) = C( u , v) \bar{u}^2 + B(u , v) \bar{v} \bar{u} + A(u , v),
\end{equation}
where for $u , v \in U$, we set
\begin{eqnarray*}
C(u , v) & = & \sum_{j=1}^m \frac{1}{|\pi_j + u^2|^2} + \frac{4 \sigma \omega%
}{|2 \omega u^2 + v u + 2 \theta|^2} > 0 , \\
B(u , v) & = & \frac{2 \sigma}{|2 \omega u^2 + v u + 2 \theta|^2} \geq 0 , \\
A(u , v) & = & \sum_{j=1}^m \frac{\pi_j}{|\pi_j + u^2|^2} + \frac{4 \sigma \theta}{|2
\omega u^2 + v u + 2 \theta|^2} > 0.
\end{eqnarray*}
We have just shown that $|R_1 (u , v)| = |2 \omega u^2 + v u + 2 \theta|>0$ if $u , v \in
U$. For $\sigma = 0$, one has
\[
R_2 (u , v) = C( u , v) \left[\bar{u}^2 + \frac{A(u , v)}{C ( u , v)} \right] \neq 0,
\]
whenever $u , v \in U$. For $\sigma \neq 0$, one has $B(u , v) > 0$ and by (%
\ref{l16}) one would get from $R_2 (u , v) = 0$
\[
\bar{v} = - \frac{A(u , v)}{B(u , v)}\cdot u - \frac{C(u , v)}{B(u , v)}%
\cdot \bar{u},
\]
which yields in turn
\[
\left(v \in U \right) \Rightarrow \left(R_2 (u , v) \neq 0 \right).
\]
\quad\hfil%
\hbox{
    \vrule width 0.8pt height 7pt depth 0.8pt  \hskip -0.2pt
    \vrule width   5pt height 7pt depth -6.2pt \hskip -5pt
    \vrule width   5pt height 0pt depth 0.8pt
    \vrule width 0.8pt height 7pt depth 0.8pt \hskip 2pt}

\textbf{Proof of Lemma \ref{l1lm}.} Similarly to (\ref{l4}) one has for $%
\varphi \in \mathcal{P}^+$,
\[
\varphi (\Delta_{\theta , \omega}) = \varphi_0 \prod_{j = 1}^m (\sigma_j + \Delta_{\theta
, \omega}), \quad \sigma_j \geq 0, \ j = 1, 2, \dots , m, \ m \in \mathbb{N}.
\]
By Lemma \ref{l2lm} each $(\sigma_j + \Delta_{\theta , \omega})$ maps $%
\mathcal{P}^+$ into itself, hence the whole $\varphi (\Delta_{\theta ,
\omega})$ does so. \quad\hfil%
\hbox{
    \vrule width 0.8pt height 7pt depth 0.8pt  \hskip -0.2pt
    \vrule width   5pt height 7pt depth -6.2pt \hskip -5pt
    \vrule width   5pt height 0pt depth 0.8pt
    \vrule width 0.8pt height 7pt depth 0.8pt \hskip 2pt}

\textbf{Proof of Theorem \ref{2.30tm}.} By (\ref{dec}) the operator $\exp(a\Delta_{\theta
, \omega })$ is a composition of $\exp(a\omega z D)$ and $\exp(\gamma \Delta_\theta )$
with $\gamma = (e^a\omega -1)/\omega$. The latter operator continuously maps
$\mathcal{L}_b^+$ in $\mathcal{L}_\beta^+$ (Proposition \ref{a1tm}), where $\beta = b[1-
( e^a\omega -1)(b/\omega)]$. The former one, also continuously, maps
$\mathcal{L}_\beta^+$ into $\mathcal{L}_{c(b)}^+$, which follows from Proposition
\ref{Pro6}.
\quad\hfil%
\hbox{
    \vrule width 0.8pt height 7pt depth 0.8pt  \hskip -0.2pt
    \vrule width   5pt height 7pt depth -6.2pt \hskip -5pt
    \vrule width   5pt height 0pt depth 0.8pt
    \vrule width 0.8pt height 7pt depth 0.8pt \hskip 2pt}

In a similar way, the proof of Theorem \ref{2.5tm} follows from (\ref{dec}), (\ref{3.1}),
and Proposition \ref{intt1pn}. The proof of Theorem \ref{2.6tm} follows from (\ref{dec}),
(\ref{3.1}), and Proposition \ref{intt2pn}. The proof of Corollary \ref{2.1co} follows
from (\ref{int5}), (\ref{int6}), and Theorem \ref{2.30tm}.

\section{Differential Equation}

Now we may use the operators introduced above to describe the solutions to certain Cauchy
problems. First we consider the following one
\begin{eqnarray}  \label{e1}
\frac{\partial f(t,z)}{\partial t} &=&(\theta +\omega z)\frac{\partial f(t,z)}{\partial
z}+z\frac{\partial ^{2}f(t,z)}{\partial z^{2}}, \ \ \ \omega \in
\mathbb{R},\ z\in \mathbb{C}, \ t>0;  \nonumber \\
f(0,z) &=&g(z) .
\end{eqnarray}

\begin{Tm}
\label{e1tm} For every $\theta \geq 0,$ $\omega \in \mathbb{R}$, and $g\in \mathcal{E}$
having the form
\begin{equation}
g(z)=\exp (-\varepsilon z)h(z),\ \ \ h\in \mathcal{A}_{0},\ \ \varepsilon \geq 0,
\label{e3}
\end{equation}
\vskip.1cm
\begin{tabular}{ll}
\textrm{(i)} & the problem (\ref{e1}) has a unique solution in $\mathcal{A}%
_{\varepsilon }$, which \\
& may be written
\end{tabular}
\begin{eqnarray}
& & f(t,z) =\left( \exp (t\Delta _{\theta, \omega })g\right) (z) =
\label{e2} \\
&=&\exp \left( -\frac{\omega z}{1-e^{-t\omega}}\right) \int_{0}^{+\infty } w_{\theta
}\left( \frac{\omega zs}{1-e^{-t\omega}}\right) g\left( s\frac{ e^{t\omega
}-1}{\omega}\right)s^{\theta -1}e^{-s} \mathrm{d}s;  \nonumber
\end{eqnarray}
\vskip.1cm
\begin{tabular}{ll}
\textrm{(ii)} & if in (\ref{e3}) $\varepsilon >0$, the solution (\ref{e2})
converges \\
& to zero in $\mathcal{A}_{\varepsilon }$ when $t\rightarrow +\infty$; \\
[0.1cm] \textrm{(iii)} & if in (\ref{e3}) $h\in \mathcal{L}_{0}\subset \mathcal{A}
_{0}$, the solution (\ref{e2}) belongs \\
& either to $\mathcal{L}_{0}$, for $\varepsilon =0$, or to $\mathcal{\ L}^{-} $, for
$\varepsilon >0$.
\end{tabular}
\end{Tm}

\textbf{Proof.} Let $\varphi _{t}(z)=\exp (tz)$. The operator valued function
$[0,t_{0})\ni t\mapsto \varphi _{t}(\Delta _{\theta, \omega})$ is continuous and
differentiable in the norm-topology, and
\[
\varphi _{t}^{\prime }(\Delta _{\theta, \omega })=\Delta _{\theta, \omega }\varphi
_{t}(\Delta _{\theta, \omega }).
\]
The functions (\ref{e3}), with a given $\varepsilon$ and all $h\in \mathcal{A}_0 $, form
a subspace of $\mathcal{B}_{b}\supset \mathcal{A}_{\varepsilon }$, $ b>\varepsilon $. The
restrictions of $\varphi _{t}(\Delta _{\theta })$, $t \in [0,t_{0})$ to this subspace is
a differentiable semi--group. Then the problem (\ref{e1}) has a unique solution in the
mentioned subspace (see e.g. Theorem 1.4 p.109 of \cite{Pazy}) having the form
\[
f(t , z) = \left(\exp(t \Delta_{\theta , \omega})g \right)(z).
\]
This proves uniqueness. The representation (\ref{e2}) follows from Theorem \ref{2.5tm}.
Further, we substitute in (\ref{e2}) the initial condition (\ref{e3}) and apply Theorem
\ref{2.6tm} with $u = -\varepsilon$. This yields
\begin{eqnarray}  \label{wz1}
f(t,z) &=&\left( 1+\varepsilon \frac{e^{t\omega}-1}{\omega}\right) ^{-\theta }\exp \left(
-\frac{\varepsilon \omega z}{\omega
e^{-t\omega }+\varepsilon (1-e^{-t\omega})}\right)  \nonumber \\
&\times &\left[ \exp \left( \frac{e^{t\omega }-1}{\omega+\varepsilon (e^{t\omega
}-1)}\Delta _{\theta }\right) h\right] \left( \frac{\omega^2 ze^{t\omega}} {\left( \omega
+\varepsilon [ e^{t\omega }-1]\right)^{2}}
\right)  \nonumber \\
&{\stackrel{\mathrm{def}}{=}}&\left( 1+\varepsilon \frac{e^{t\omega}-1}{\omega
}\right)^{-\theta }\exp \left( -\frac{\varepsilon \omega z}{\omega
e^{-t\omega}+\varepsilon (1-e^{-t\omega })}\right) h_{t}\left( z\right) . \nonumber
\end{eqnarray}
Since $h \in \mathcal{A}_0$, by Theorem \ref{2.6tm}, $h_{t}\in \mathcal{A}_{0}$, and by
Corollary \ref{2.1co}, $h_{t}\in \mathcal{L}_{0}$ if $h\in \mathcal{L}_{0}$. The former
yields that the solution belongs to $\mathcal{A}_{\varepsilon }$ and the latter does
claim (iii). It remains to prove the convergence stated in (ii).  The continuity of the
operator $\exp (t\Delta _{\theta, \omega })$ yields that in $\mathcal{A}_{0}$
\begin{eqnarray*}
h_{t}(z) &=&\left[ \exp \left( \frac{e^{t\omega}-1}{b+\varepsilon (e^{t\omega}-1)} \Delta
_{\theta }\right) h\right] \left( \frac{ze^{t\omega }
}{\left( 1+\varepsilon\omega^{-1} [e^{t\omega}-1]\right) ^{2}}\right) \\
\ &\rightarrow &\ \left\{ \exp \left( \frac{1}{\varepsilon }\Delta _{\theta }\right)
h\right\} (0),
\end{eqnarray*}
if $\omega \neq 0$. The case $\omega =0$ may be handled similarly. Therefore, the product
in (\ref{wz1}) tends to zero in $\mathcal{A}_{\varepsilon }$ when $t\rightarrow +\infty $.
\quad\hfil%
\hbox{
    \vrule width 0.8pt height 7pt depth 0.8pt  \hskip -0.2pt
    \vrule width   5pt height 7pt depth -6.2pt \hskip -5pt
    \vrule width   5pt height 0pt depth 0.8pt
    \vrule width 0.8pt height 7pt depth 0.8pt \hskip 2pt}

Given $N \in \mathbb{N}$, let $\mathcal{E}^{(N)}$ stand for the set of analytic functions
$F: \mathbb{R}^{N}\rightarrow \mathbb{C}$. For $b>0$, we set
\begin{equation}
\Vert F\Vert _{b,N}\ \stackrel{\mathrm{def}}{=}\ \sup_{x\in \mathbb{R}^{N}} \{\mid
F(x)\mid \exp (-b\mid x\mid ^{2})\},  \label{9}
\end{equation}
where $\mid x\mid $ is the Euclidean norm. Set
\begin{equation}
\mathcal{A}_{a}^{(N)}\ \stackrel{\mathrm{def}}{=}\ \{F\in \mathcal{E}^{(N)}\ \mid \ \Vert
F\Vert _{b,N}<\infty ,\ \forall b>a\},\ \ a\geq 0.  \label{10}
\end{equation}
This set equipped with the topology generated by the family of norms $%
\{\Vert .\Vert_{b,N},b>a\}$ becomes a Fr{\'{e}}chet space. Let $O(N)$ stand for the group
of all orthogonal transformations of $\mathbb{R}^{N}$. A function $F\in
\mathcal{E}^{(N)}$ is said to be isotropic if for every $U\in
O(N)$ and all $x\in \mathbb{R}^{N}$, one has $F(Ux)=F(x)$. The subset of $%
\mathcal{E}^{(N)}$ consisting of isotropic functions is denoted by $\mathcal{%
E}_{\mathrm{isot}}^{(N)}$. Let also $\mathcal{P}_{\mathrm{isot}%
}^{(N)}\subset \mathcal{E}_{\mathrm{isot}}^{(N)}$ stand for the set of isotropic
polynomials. The classical Study--Weyl theorem (see \cite{Lu}) implies that there exists
a bijection between the set of all polynomials of a single complex variable $\mathcal{P}$
and $\mathcal{P}_{ \mathrm{isot}}^{(N)}$ established by
\[
\mathcal{P}_{\mathrm{isot}}^{(N)}\ni P(x)=p((x,x))\in \mathcal{P},
\]
where $(.,.)$ is the scalar product in $\mathbb{R}^{N}$. Obviously, each a function $F$
having the form
\begin{equation}
F(x)=f((x,x)),  \label{11x}
\end{equation}
with a certain $f\in \mathcal{E}$, belongs to $\mathcal{E}_{\mathrm{isot}
}^{(N)}$. Given $\mathcal{X}\subset \mathcal{E}$, we write $%
\mathcal{X}(\mathbb{R}^{N})$ for the subset of $\mathcal{E}_{\mathrm{isot}
}^{(N)}$ consisting of the functions obeying (\ref{11x}) with $f\in \mathcal{%
\ X}$. Consider
\[
\mathcal{E}_{\mathrm{isot}}^{(N)}\ni F\mapsto \left( \Delta +\left( \frac{d}{%
(x,x)}+b\right) (x,\nabla )\right) F\in \mathcal{E}_{\mathrm{isot}}^{(N)},
\]
where $\Delta $ and $\nabla $ are the Laplacian and gradient in $\mathbb{R}^{N}$. For $F$
and $f$ satisfying (\ref{11x}), one has
\begin{eqnarray}  \label{17}
\left( \Delta +\left( \frac{d}{(x,x)}+b\right) (x,\nabla )\right) F(x)  =  4\left(
\Delta_{\theta , \omega} f\right) ((x , x)),
\end{eqnarray}
where $\Delta _{\theta }$ and $\Delta_{\theta, \omega}$ are defined by (\ref{1.1}) with
\begin{equation}
\theta =\frac{N+d}{2}, \qquad \omega =\frac{b}{4}.  \label{16}
\end{equation}
Consider the following Cauchy problem
\begin{eqnarray}
\frac{\partial F(t,x)}{\partial t} &=&\left( \Delta +\left( \frac{d}{(x,x)}%
+b\right) (x,\nabla )\right) F(t,x),  \label{luna1} \\
F(0,x) &=&G(x)\in \mathcal{E}_{\mathrm{isot}}^{(N)}.  \nonumber
\end{eqnarray}
where $t\in \mathbb{R}_{+}\ $and $x\in \mathbb{R}^{N}.$

\begin{Tm}
\label{Lunatm}For every $d\geq -N,$ $b\in \mathbb{R}$, and $G$ having the form
\begin{equation}
G(x)=\exp [-\varepsilon (x,x)]h((x,x)),\ \ h\in \mathcal{A}_{0},\ \ \varepsilon \geq 0,
\label{luna2}
\end{equation}

\begin{tabular}{ll}
\textrm{(i)} & the problem (\ref{luna1}) has a unique solution in $\mathcal{A%
}_{\varepsilon }^{(N)}$, which, \\
& for $t>0$, may be written as follows,
\end{tabular}
\begin{eqnarray}  \label{luna3}
F(t,x) &=&\exp \left( -\frac{b(x,x)}{4(1-e^{-tb})}\right) \int_{0}^{+\infty
}w_{\theta }\left( \frac{bs(x,x)}{4(1-e^{-tb})}\right) \\
& & h\left(4s\frac{e^{tb}-1}{b}\right)s^{\theta -1}\exp \left(-s\left(1+4\varepsilon
\frac{e^{tb}-1}{b}\right)\right)\mathrm{d}s, \nonumber
\end{eqnarray}
where $\theta $ is given by (\ref{16});

\begin{tabular}{ll}
\textrm{(ii)} & if in (\ref{luna2}) $\varepsilon >0$, the solution (\ref
{luna3}) converges to zero in $\mathcal{A}_{\varepsilon }^{(N)}$ \\
& when $t\rightarrow +\infty $; \\[0.1cm]
\textrm{(iii)} & if in (\ref{luna2}) $h\in \mathcal{L}_{0}\subset \mathcal{A}%
_{0}$, the solution (\ref{luna3}) \\
& belongs either to $\mathcal{L}_{0}(\mathbb{R}^{N})$, for $\varepsilon =0$, or to
$\mathcal{L}^{-}(\mathbb{R}^{N})$, for $\varepsilon >0$.
\end{tabular}
\end{Tm}

The proof directly follows from Theorem \ref{e1tm} on the base of the
correspondence formulas (\ref{11x}) and (\ref{17}).


\begin{thebibliography}{00}

\bibitem{Il}  L. Iliev, "Laguerre Entire Functions", Bulgarian Academy of
Sciences, Sofia, 1987.

\bibitem{KOU}  Yu. Kozitsky, P. Oleszczuk and G. Us, Integral operators and
dual orthogonal systems on a half-line, \textit{Integral Transform. Spec. Funct.}, 12
(2001), 257-278.

\bibitem{Kor} Ju. F. Korobeinik, On the question of the representation of an arbitrary
linear operator in the form of a differentail operator of infinite order, \textit{Mat.
Zametki} 16 (1974), 277-283.

\bibitem{KozW1}  Yu. Kozitsky, L. Wo{\l}owski, Laguerre entire functions and
related locally convex spaces, \textit{Complex Variables Theory Appl.,} 44 (2001),
225-244.

\bibitem{Lev}  B. J. Levin, "Distribution of Zeros of Entire Functions",
Amer. Math. Soc., 1964.

\bibitem{Lu}  D. Luna, Fonctions diff{\'e}rentiables invariantes sous l'op{%
\'e}rations d'un groupe r{\'e}ductif, \textit{Ann. Inst. Fourier, Grenoble,} 26 (1976),
33-49.

\bibitem{Nzag}  H. Neidthard, V.A. Zagrebnov, Trotter-Kato product formula
and operator-norm convergence, \textit{Comm. Math. Phys.,} 205 (1999), 129-159.

\bibitem{Pazy}  A. Pazy, "Semi-Groups of Linear Operators and Applications
to Partial Differential Equations", Lecture Notes of the University of Maryland, College
Park, 1974.

\end{thebibliography}
\end{document}